\documentclass[reqno]{amsart}

\usepackage{amssymb,latexsym,mathrsfs,amsmath}
\usepackage{amsthm}
\usepackage{graphicx}
\usepackage{diagbox}
\usepackage{mathrsfs}
\usepackage{tikz}
\usetikzlibrary{arrows}

\newtheorem{lemma}{Lemma}[section]

\newtheorem{theorem}[lemma]{Theorem}

\theoremstyle{definition}
\newtheorem{remark}[lemma]{Remark}
\newtheorem{definition}[lemma]{Definition}

\DeclareMathOperator{\id}{id}

\DeclareMathOperator{\Lee}{Lee}
\DeclareMathOperator{\BN}{BN}
\DeclareMathOperator{\xo}{xo}

\DeclareMathOperator{\tpg}{\widetilde{pg}}

\newcommand{\Z}{\mathbb{Z}}

\newcommand{\Q}{\mathbb{Q}}

\newcommand{\F}{\mathbb{F}}

\newcommand{\Th}{T^{\frac{1}{2}}}

\newcommand{\CL}{C_{\Lee}}
\newcommand{\HL}{H_{\Lee}}
\newcommand{\CBN}{C_{\BN}}
\newcommand{\HBN}{H_{\BN}}
\newcommand{\tCBN}{\widetilde{C}_{\BN}}
\newcommand{\tHBN}{\widetilde{H}_{\BN}}
\newcommand{\FrobL}{\mathcal{F}_{\Lee}}
\newcommand{\FrobCL}{\overline{\mathcal{F}}_{\Lee}}
\newcommand{\FrobB}{\mathcal{F}_{\BN}}

\begin{document}
\parindent0em
\setlength\parskip{.1cm}
\thispagestyle{empty}
\title{A note on the $X$-torsion order of a knot}
\author[Dirk Sch\"utz]{Dirk Sch\"utz}
\address{Department of Mathematical Sciences\\ Durham University\\ Durham DH1 3LE\\ United Kingdom}
\email{dirk.schuetz@durham.ac.uk}
%\subjclass[2020]{primary: 57K18; secondary: 57K10}
%\keywords{Odd Khovanov homology, Steenrod square}

\begin {abstract}
We show that the $X$-torsion order of a knot, which is defined in terms of a generalised Lee complex, can be calculated using the reduced Bar-Natan--Lee--Turner spectral sequence. We use this for extensive calculations, including an example of $X$-torsion order $4$.
\end{abstract}

\maketitle

\section{Introduction}
In recent years torsion invariants for knots arising from Khovanov homology appeared in many contexts. For example, they give rise to lower bounds for unknotting numbers and Gordian distance \cite{MR4007369, MR3894208, lewark2024khov}, rational unknotting numbers \cite{iltgen2021khovanov}, band unknotting numbers \cite{MR4420588}, ribbon distance \cite{MR4092319, gujral2020ribbon}, and Turaev genus \cite{MR4417719}.

They are defined from deformations of Khovanov homology over a polynomial ring with coefficients in a field $\F$, using the maximal torsion order of the homology viewed as a module over the polynomial ring. These deformations come with spectral sequences that start with the Khovanov homology of the knot and converge to the (shifted) Khovanov homology of the unknot. There is a close relation between the torsion order and the number of pages of the spectral sequence. The odd one out here is the $X$-torsion order $\xo_\F(K)$, which only satisfies $\xo_\F(K)\in \{2k-3, 2k-2\}$ with $k$ the number of pages in the Lee spectral sequence (we assume that the $E_1$-page is the Khovanov homology, and $E_k$ the first page equal to $E_\infty$). By using the $X$-action on the Lee complex, one can get a spectral sequence where this is more precise, but our main theorem states that we can also use the reduced Bar-Natan--Lee--Turner spectral sequence.

\begin{theorem}\label{thm:main}
Let $K$ be a knot and $\F$ a field of characteristic different from $2$. Then
\[
\xo_\F(K) = \tpg_\F(K)-1,
\]
where $\tpg_\F(K)$ is the number of pages in the reduced Bar-Natan--Lee--Turner spectral sequence of $K$ with coefficients in $\F$.
\end{theorem} 

%This also works for links, but in this case both the $X$-torsion order and the reduced spectral sequence require the choice of a base point on the link.

{\bf Acknowledgements:} The author would like to thank Nathan Dunfield for valuable discussions on the $X$-torsion order.

\section{Lee and Bar-Natan homology}
A (commutative) \em Frobenius system \em is a tuple $\mathcal{F}=(R,A,\varepsilon, \Delta)$ with $A$ a commutative ring and a subring $R$, $\varepsilon\colon A\to R$ an $R$-module map, $\Delta\colon A\to A\otimes_R A$ an $A$-bimodule map that is co-associative and co-commutative, such that $(\varepsilon\otimes \id)\circ \Delta = \id$.

Given a Frobenius system $\mathcal{F}= (R,A,\varepsilon,\Delta)$ such that $A$ is free of rank $2$ over $R$, Khovanov  \cite{MR2232858} showed that for a link diagram $D$ one can define a cochain complex $C(D;\mathcal{F})$ over $R$ whose homology is a link invariant.

The classical example of Khovanov homology is obtained by choosing $R=\Z$, $A=Z[X]/(X^2)$, $\varepsilon(1) = 0$ and $\varepsilon(X) = 1$, and $\Delta$ is given by $\Delta(1) = 1\otimes X + X\otimes 1$. With $\Delta(X) = X\otimes X$ one has indeed that $\Delta$ is an $A$-bimodule map.

We will be mainly interested in two deformations of this system, which we name after Lee and Bar-Natan.

\subsection{Lee homology}
The Lee deformation of Khovanov homology can be described as follows. For the ground ring we use $R[T]$, and $A = R[X,T]/(X^2-T)$. Both rings are graded by declaring $\deg(1) = 0$, $\deg(T) = -4$ and $\deg(X) = -2$. The Frobenius system is given by
\[
\Delta_A(1) = 1\otimes X + X\otimes 1 \hspace{2cm} \Delta_A(X) = X\otimes X+ T\otimes 1
\]
and co-unit $\varepsilon\colon A \to R[T]$ given by $\varepsilon(1) = 0$, $\varepsilon(X) = 1$. Note that since tensor products are over $R[T]$, $T\otimes 1 = 1\otimes T$.

We denote this Frobenius system by $\FrobL$ and the resulting link homology chain complex for a link diagram $D$ by $\CL(D;R[T])$. Using creative grading shifts, one can ensure that this complex is bigraded, that is, it has a grading different from the homological grading, called the $q$-grading, which is preserved by the boundary.

As an $R[T]$-module, $A$ is free of rank $2$, with a basis given by $\{1,X\}$. Since $T= X^2$ in $A$, we have $A\cong R[X]$. By choosing a base point on $D$, we can turn $\CL(D;R[T])$ into an $R[X]$-chain complex.

If $R=\F$ is a field, then $\F[X]$ is a Euclidean domain, and therefore
\begin{equation}\label{eq:SmithNormalLee}
\CL(D;\F[T]) \cong \bigoplus_{i\in I} D_i,
\end{equation}
where $I$ is a finite set, $D_i$ is either a single copy of $\F[X]$, or $D_i$ is concentrated in two adjacent homological degrees, and is of the form
\[
\F[X] \stackrel{X^{k_i}}{\longrightarrow} \F[X],
\]
where $k_i$ is a non-negative integer. This follows from the usual Smith-Normalization process, noting that we can keep this grading-preserving at every step. Also, we allow $k_i = 0$, so that we get an isomorphism of $\F[X]$ complexes in (\ref{eq:SmithNormalLee}). As $T=X^2$, we can also view this as an isomorphism of $\F[T]$ complexes.

Note that the free part for a knot is just one copy, see Turner \cite{MR4079621}. More generally, if $L$ is a $c$-component link and $\F$ a field of characteristic different from $2$, the homology of $\CL(D;\F[T])$ decomposes into $\F[X]^{2^{c-1}}\oplus T(L)$, where
\[
T(L) = \{ a\in \HL(L;\F[T])\mid X^n a = 0 \mbox{ for some }n\}.
\]

\begin{definition}
Let $L$ be a link with basepoint and $\F$ a field of characteristic different from $2$. The {\em $X$-torsion order of $L$}, denoted $\xo_\F(L)$, is defined as the minimal $n$ such that $X^nT(L) = 0$.
\end{definition}

Notice that $\xo_\F(L)$ is the largest $k_i$ that appears in the decomposition (\ref{eq:SmithNormalLee}).

\begin{remark}
For a knot the basepoint does not affect $\xo_\F(L)$, but for a link it can make a difference. The easiest way to see this is to consider a split link with an unknot component.
\end{remark}

Lee \cite{MR2173845} originally worked over $\Q$. To get the cochain complex from \cite{MR2173845}, we only need to use the change of base ring homomorphism $\eta\colon \Q[T]\to\Q$ sending $T$ to $1$.

\subsection{Bar-Natan homology}
To get the Bar-Natan deformation of Khovanov homology we write the ground ring as $R[H]$ and use $B = R[X,H]/(X^2-XH)$, with co-multiplication given by
\[
\Delta_B(1) = 1\otimes X +X\otimes 1 - H\otimes 1 \hspace{2cm} \Delta_B(X) = X\otimes X,
\]
and co-unit $\varepsilon\colon B\to R[H]$ given by $\varepsilon(1) = 0$, $\varepsilon(X) = 1$.

We write $\FrobB$ for this Frobenius system. Again we get gradings on the ground ring and $B$ by setting $\deg(1) = 0$, $\deg(H) = -2 = \deg(X)$.

For a link diagram $D$ we denote the resulting chain complex by $\CBN(D; R[H])$ and the homology by $\HBN(D;R[H])$. Again this is bigraded.

As in the case of the Lee complex, when viewed as an $R[H]$-module, $B$ is free of rank $2$ with basis given by $\{1,X\}$. However, as $X (X-H) = 0$ in $B$, we do not get an isomorphism with $R[X]$. 

The advantage of the Bar-Natan complex is that it behaves better over $\Z$, and in particular, over $\F_2$, the field with two elements. If $\F$ is a field of characteristic different from $2$, it is also closely related to the Lee complex, as we will see in Section \ref{sec:mainproof}.

To explain what we mean by `behaving better' than the Lee complex, consider the following two observations.
\begin{itemize}
\item If $\eta\colon \Z[H]\to S$ is a ring homomorphism such that $\eta(H)$ is a unit, then $\HBN(K;S) = H(\CBN\otimes_{\Z[H]}S)\cong S\oplus S$, concentrated in homological degree $0$. This follows from \cite[Prop.2.1]{MR4504654}.
\item There is a well defined reduced complex
\[
\tCBN(D;\Z[H]) = X\cdot \CBN(D;\Z[H])
\]
after choosing a basepoint on the link diagram $D$. With $\eta$ as above we get $\tHBN(L;S) = H(\tCBN\otimes_{\Z[H]} S) \cong S^{2^{c-1}}$ for a $c$-component link $L$.
\end{itemize}

\section{Proof of the Main Theorem}
\label{sec:mainproof}

In this section, $\F$ is a field of characteristic different from $2$.

Let us introduce a formal variable $\Th$ with $(\Th)^2 = T$ and consider the inclusion $\F[T]\to \F[\Th]$. We get a new Frobenius system $\FrobCL = \FrobL\otimes_{\F[T]}\F[\Th]$, that is, the ground ring is $\F[\Th]$, $A$ is given by $\F[X,\Th]/(X^2-T)$, and $\varepsilon$ and $\Delta$ are as in $\FrobL$. The resulting link complex is denoted by $\CL(D;\F[\Th]) = \CL(D;\F[T])\otimes_{\F[T]}\F[\Th]$.

We can define a ring isomorphism
\[
\Phi\colon \F[X,H]/(X^2-XH) \to \F[X,\Th]/(X^2-T)
\]
by sending $H$ to $\Th$ and $X$ to $\frac{1}{2}(X+\Th)$. Note that this sends $X-H$ to $\frac{1}{2}(X-\Th)$ (making it a well defined ring homomorphism) and $2X-H$ to $X$.

This does not quite induce an isomorphism of Frobenius systems $\FrobB$ to $\FrobCL$, but it does if we twist the latter by $\frac{1}{2}$, compare \cite{MR2232858}. In particular, it induces a grading preserving isomorphism of $\F[H]$-cochain complexes
\[
\Phi\colon \CBN(D;\F[H]) \to \CL(D;\F[\Th]).
\]
Here we treat the latter as an $\F[H]$-complex using the identification of $\F[H]$ and $\F[\Th]$ by restricting $\Phi$ to $\F[H]$.

Recall the isomorphism (\ref{eq:SmithNormalLee}) and treat it as an isomorphism over $\F[T]$. It induces an isomorphism of $\F[\Th]$-cochain complexes
\[
\CL(D;\F[\Th])\cong \bigoplus_{i\in I} D_i\otimes_{\F[T]}\F[\Th].
\]
Combining with $\Phi$ and using that $\Phi(2X-H) = X$ this shows that $\CBN(D;\F[H])$ to a direct summand of complexes $\overline{D}_i$, each either being a single free copy $\F[X,H]/(X^2-XH)$, or a complex of the form
\[
\begin{tikzpicture}
\node at (0,0) {$\F[X,H]/(X^2-XH)$};
\node at (6,0)  {$\F[X,H]/(X^2-XH)$};
\draw[->] (1.8,0) --node [above] {$(2X-H)^{k_i}$} (4.2,0);
\end{tikzpicture}
\]
with $k_i\geq 0$. Passing to the reduced complex shows that $\tCBN(D;\F[H])$ is isomorphic to a direct sum of free copies of $\F[H]$ and complexes of the form
\[
\begin{tikzpicture}
\node at (0,0) {$\F[H]$};
\node at (2.4,0)  {$\F[H].$};
\draw[->] (0.6,0) --node [above] {$H^{k_i}$} (1.8,0);
\end{tikzpicture}
\]
From this isomorphism it follows that the spectral sequence starting with the reduced Khovanov homology with coefficients in $\F$ collapses after $k$ steps, where $k$ is the maximum of the $k_i$. Since this maximum is also $\xo_\F(K)$, as follows from (\ref{eq:SmithNormalLee}), this proves Theorem \ref{thm:main}.

\begin{remark}
Viewing the Lee complex $\CL(D;\F[T])$ of a link with base point as an $\F[X]$-complex and using $\eta\colon \F[X]\to \F$ with $\eta(X)=1$, gives rise to a filtered complex $\CL(D;\F[T])\otimes_{\F[X]}\F$. Our argument shows that the corresponding spectral sequence agrees with the reduced Bar-Natan--Lee--Turner spectral sequence.
\end{remark}

\section{Computations}

In view of Theorem \ref{thm:main} we extend the definition of $X$-torsion for $\F_2$ as follows.

\begin{definition}
Let $L$ be a link with a basepoint. We then define
\[
\xo_{\F_2}(L) = \tpg_{\F_2}(L)-1,
\]
where $\tpg_{\F_2}(L)$ is the number of pages in the reduced Bar-Natan--Lee--Turner spectral sequence of $K$ with coefficients in $\F_2$.
\end{definition}

This agrees with the $H$-torsion order of \cite{MR4007369}.

Computations of $\tpg_{\F}(K)$ are readily available, for example, using \verb+knotjob+, which can be found at the author's website.

For knots with up to $14$ crossings, $\xo_\F(K)$ does not depend on $\F$. But there are five knots with $15$ crossings where $\xo_\Q(K) > \xo_{\F_2}(K)$ and one $15$-crossing knot with $\xo_\Q(K) < \xo_{\F_2}(K)$. There are also $111$ knots with $16$ crossings such that $\xo_\Q(K)\!= \xo_{\F_2}(K)$. For all knots with up to 16 crossings we have $\xo_\Q(K) = \xo_{\F_3}(K)$, and no $X$-torsion order is bigger than $2$.

The Manolescu--Marengon knot $K$, which is a counterexample to the Knight-move conjecture \cite{MR4042864}, satisfies $\xo_\Q(K) = 3$. Interestingly, we get $\xo_{\F_2}(K) = \xo_{\F_3}(K) = 2$ for this knot. In particular, this knot satisfies the Knight-move conjecture in characteristics $2$ and $3$. There is a slight variation of this knot, $K'$, which uses a full twist on $8$ strands as opposed to $6$ strands. Calculations show this knot satisfies $\xo_\Q(K') = 4$. The next variation, which would use a full twist on $10$ strands, is unfortunately outside of the range for calculations.

For torus knots examples with $X$-torsion order bigger than $2$ for finite fields are known in the cases $T(5,6)$ and $T(7,8)$, compare \cite{MR2320156, MR4417719}. More can be said about torus knots, but we only want to highlight the case $T(8,9)$, where $\xo_\Q(T(8,9)) =2 < 3 = \xo_{\F_7}(T(8,9))$, despite the Betti numbers of Khovanov homology being the same with $\Q$ and $\F_7$ coefficients.

\bibliography{KnotHomology}
\bibliographystyle{amsalpha}

\end{document}